\documentclass{amsart}
\usepackage{amscd,amssymb,stmaryrd}
\usepackage[all]{xy}
\renewcommand{\mod}{\operatorname{mod}\nolimits}
\newcommand{\umod}{\operatorname{\underline{mod}}\nolimits}
\newcommand{\gr}{{\operatorname{gr}\nolimits}}
\newcommand{\add}{\operatorname{add}\nolimits}
\newcommand{\Hom}{\operatorname{Hom}\nolimits}
\newcommand{\uHom}{\operatorname{\underline{Hom}}\nolimits}
\renewcommand{\Im}{\operatorname{Im}\nolimits}
\newcommand{\Ker}{\operatorname{Ker}\nolimits}

\newcommand{\rrad}{\mathfrak{r}}
\newcommand{\Ann}{\operatorname{Ann}\nolimits}
\newcommand{\Tr}{\operatorname{Tr}\nolimits}
\newcommand{\Ext}{\operatorname{Ext}\nolimits}
\newcommand{\cExt}{\operatorname{\widehat{Ext}}\nolimits}
\newcommand{\op}{{\operatorname{op}\nolimits}}

\newcommand{\id}{{\operatorname{id}\nolimits}}

\newcommand{\m}{\mathfrak{m}}
\newcommand{\frakp}{\mathfrak{p}}
\newcommand{\fraka}{\mathfrak{a}}

\newcommand{\G}{\Gamma}
\renewcommand{\L}{\Lambda}
\newcommand{\Z}{{\mathbb Z}}
\newcommand{\C}{{\mathcal C}}
\newcommand{\E}{{\mathcal E}}
\newcommand{\N}{{\mathcal N}}

\newcommand{\Y}{{\mathcal Y}}
\renewcommand{\P}{{\mathcal P}}

\newcommand{\extto}{\xrightarrow}
\newcommand{\MaxSpec}{\operatorname{MaxSpec}\nolimits}
\newcommand{\HH}{\operatorname{HH}\nolimits}

\newcommand{\Ind}{\operatorname{Ind}\nolimits}
\newtheorem{lem}{Lemma}[section]
\newtheorem{prop}[lem]{Proposition}
\newtheorem{cor}[lem]{Corollary}
\newtheorem{thm}[lem]{Theorem}
\theoremstyle{definition}
\newtheorem{defin}[lem]{Definition}
\newtheorem*{remark}{Remark}
\newtheorem{assumption}{Assumption}

\begin{document}

\title{Support varieties for selfinjective algebras}
\author[Erdmann]{Karin Erdmann}
\address{Karin Erdmann, Miles Holloway, Rachel Taillefer\\ 
Mathematical Institute\\ 24--29 St.\ Giles\\
Oxford OX1 3LB\\ England}
\email{erdmann, holloway \textrm{and} taillefe \textrm{at} maths.ox.ac.uk}
\author[Holloway]{Miles Holloway}
\author[Snashall]{Nicole Snashall}
\address{Nicole Snashall\\ Department of Mathematics\\
University of Leicester\\ University Road\\ Leicester, LE1 7RH\\ England}
\email{njs5@mcs.le.ac.uk}
\author[Solberg]{\O yvind Solberg}
\address{\O yvind Solberg\\Institutt for matematiske fag\\
NTNU\\ N--7491 Trondheim\\ Norway}
\email{oyvinso@math.ntnu.no}
\author[Taillefer]{Rachel Taillefer}
\begin{abstract}
Support varieties for any finite dimensional algebra over a field were 
introduced in \cite{SS} using graded subalgebras of the Hochschild
cohomology. We mainly study these varieties for selfinjective algebras
under appropriate finite generation hypotheses. Then many of the
standard results from the theory of support varieties for finite
groups generalize
to this situation. In particular, the complexity of the module equals 
the dimension of its corresponding variety, all closed homogeneous
varieties occur as the variety of some module, the variety of an
indecomposable module is connected, periodic modules are lines and for 
symmetric algebras a generalization of Webb's theorem is true. 
\end{abstract}
\date{\today}
\maketitle

\section*{Introduction}

Let $k$ be a field of characteristic $p$ and $G$ a finite group.  In
1971, Quillen \cite{Q} gave a description of the cohomology ring
$H^*(G,k)$ modulo nilpotent elements as an inverse limit of cohomology
algebras of elementary abelian $p$-subgroups of $G$. This has led to
the work of Benson, Carlson and others on the theory of varieties for
$kG$-modules, and in general to deep structural information about
modular representations of finite groups \cite{B}, \cite{C}.
Subsequently, analogous results have been obtained for $p$-Lie
algebras (by Friedlander and Parshall \cite{FP}, Jantzen \cite{J}, and
others), and also for Steenrod algebras arising in algebraic topology
(by Palmieri \cite{Pa}).

The support variety of a $kG$-module is a powerful invariant.  This is
defined in terms of the maximal ideal spectrum of the group cohomology
$H^*(G,k)$, a finitely generated (almost) commutative graded ring.  It
acts on ${\rm Ext}^*(M,M)$ for any finitely generated $kG$-module $M$
and the support variety of the module is the variety associated to the
annihilator ideal of this action.  This construction is based on the
Hopf algebra structure which is generally not available.  In \cite{SS}
an analogous construction for an arbitrary finite dimensional algebra
$\L$ was developed where instead of group cohomology, one takes the
Hochschild cohomology ring $\HH^*(\L)$.

In this paper, we mainly study these varieties for selfinjective
algebras.  We prove that under appropriate finite generation
hypotheses many of the properties known for the group algebra
situation have analogues.

Throughout this paper $\L$ always denotes an indecomposable finite
dimensional algebra over an algebraically closed field $k$, with
Jacobson radical $\rrad$.  Recall from \cite{SS} that the variety of a
finitely generated left $\L$-module $M$ relative to a (Noetherian)
graded subalgebra $H$ of $\HH^*(\L)$ is given by
\[ V_H(M)=\{\m\in \MaxSpec H\mid \Ann_H\Ext_\L^*(M,\L/\rrad)\subseteq
\m\},\] where $\MaxSpec H$ is the maximal ideal spectrum of $H$ and
$\Ann_H X$ is the annihilator of an $H$-module $X$. In the theory of
support varieties for group rings of finite groups (\cite{B,C2,C}),
for more general finite dimensional cocommutative Hopf algebras
(\cite{FS}), and for complete intersections (\cite{Av,AB}), the
property of having a Noetherian ring $H$ of cohomological operators
over which the extension groups $\Ext^*(M,N)$ are finitely generated
$H$-modules for all finitely generated modules $M$ and $N$, is one of
the corner stones for the whole theory. Hence, two assumptions are
central: \textbf{Fg1}: $H$ is a commutative Noetherian graded
subalgebra of $\HH^*(\L)$ with $H^0=\HH^0(\L)$, and \textbf{Fg2}:
$\Ext^*_\L(M,N)$ is a finitely generated $H$-module for all finitely
generated left $\L$-modules $M$ and $N$ (see section \ref{section1}). 

In the first section we analyse the consequences of these assumptions,
and show that the algebra $\L$ must be Gorenstein (that is, the injective
dimensions of $\L$ as a left and as a right module are finite).
Furthermore, the dimension of the variety of a module is given by the
complexity of the module, and the variety of a module is trivial if
and only if the module has finite projective dimension.

The second section is devoted to characterizing elements in the
annihilator of $\Ext^*_\L(M,M)$ as an $\HH^*(\L)$-module. In the
process we introduce for each homogeneous element $\eta$ in
$\HH^*(\L)$ a bimodule $M_\eta$ that we use in the next section with
assumptions \textbf{Fg1} and \textbf{Fg2} to show that any closed
homogeneous variety occurs as the variety of some module (Theorem
\ref{thm:allhomvar}). To our knowledge the proof we give of this fact
also gives an alternative proof of the same result for support
varieties for group rings of finite groups. The proof in the group
ring case uses rank varieties and restriction to elementary abelian
subgroups. We do not yet have an analogue of a rank variety to offer
in our more general setting, but there are partial answers in special
cases considered by Erdmann and Holloway in \cite{EH}.

Again with assumptions \textbf{Fg1} and \textbf{Fg2}, periodic modules
for selfinjective algebras are characterized in the fourth section as
modules with complexity one, as in the group ring case. Using this we
prove a generalisation of Webb's theorem to selfinjective algebras
where the Nakayama functor is of finite order for any indecomposable
module.

The next section is devoted to briefly discussing relationships
between representation type and complexity. In the final section we
show that the variety of an indecomposable module is connected
whenever \textbf{Fg1} and \textbf{Fg2} hold.

We end this introduction by setting the notation and the overall
general assumptions. For any ring $R$ we denote by $\mod R$ the
finitely presented left $R$-modules. Recall that throughout the paper
$\L$ denotes an indecomposable finite dimensional algebra over an
algebraically closed field $k$ with Jacobson radical $\rrad$.  The
Hochschild cohomology ring of $\L$ decomposes in the same way as $\L$
decomposes as an algebra, so we can without loss of generality assume
throughout that $\L$ is an indecomposable (connected) algebra.  The
stable category of $\mod\L$ is denoted by $\umod\L$, that is, $\mod\L$
modulo the ideal given by the morphisms factoring through projective
modules in $\mod\L$.  We denote the enveloping algebra $\L\otimes_k
\L^\op$ by $\L^e$, and we view $\L$-bimodules as left $\L^e$-modules.
For a module $X$ in $\mod\L$ the full subcategory $\add X$ is given by
all the direct summands of every finite direct sum of copies of $X$.
For a homogeneous ideal $\mathfrak{a}$ in $H$ let $V_H(\mathfrak{a})$
be the subvariety given by the maximal ideals of $H$ containing
$\mathfrak{a}$. Let $D=\Hom_k(-,k)\colon \mod\L\to \mod\L^\op$ be the
usual duality. For a $\L$-module $B$ denote by $B^*$ the right
$\L$-module $\Hom_\L(B,\L)$.

\section*{Acknowledgements}

The second author would like to thank EPSRC for their support through
a Postdoctoral Fellowship while this research was carried out.  The
third and the fourth authors would like to express their gratitude for
the hospitality and support from the Mathematical Institute,
University of Oxford. The last author kindly acknowledges financial
support from HEFCE through a Career Development Fellowship.

\section{Finite generation}\label{section1}
Given two $\L$-modules $M$ and $N$ in $\mod\L$ the direct sum
$\Ext^*_\L(M,N)$ is a left and a right module over the Hochschild
cohomology ring $\HH^*(\L)$ of $\L$, where the left and the right
actions are related in a graded commutative way (see \cite[Theorem
1.1]{SS}).  Furthermore, it is known that in general $\HH^*(\L)$ is
not a finitely generated algebra over $k$ and that $\Ext^*_\L(M,N)$ is
not a finitely generated module over $\HH^*(\L)$. This section is
devoted to investigating consequences of finite generation of
$\Ext^*_\L(M,N)$ as a module over commutative Noetherian graded
subalgebras $H$ of $\HH^*(\L)$ for specific pairs and for all pairs of
modules in $\mod\L$. In particular, we show that if
$\Ext^*_\L(D(\L_\L),\L/\rrad)$ and $\Ext^*_\L(\L/\rrad,\L)$ are
finitely generated modules over such an $H$, then $\L$ is a Gorenstein
algebra.

As we have pointed out already, even seemingly weak finite generation
assumptions imply strong conditions on the algebras we consider. To
see how severely finite generation can fail in general, we pause to
consider the following example. Let $\L=
k\langle\alpha_1,\ldots,\alpha_n\rangle/
(\{\alpha_i\alpha_j\}_{i,j=1}^{n,n})$ for some $n\geq 2$. This is a
finite dimensional Koszul algebra, and the Koszul dual is given by
$E(\L)=\Ext^*_\L(\L/\rrad,\L/\rrad)\simeq k\langle
\alpha_1^*,\ldots,\alpha_n^*\rangle$, a free algebra in $n$
indeterminates. In \cite{SS} it is shown that the image of the map
$-\otimes_\L \L/\rrad\colon \HH^*(\L)\to E(\L)$ is contained in the
graded centre $Z_\gr(E(\L))$ of $E(\L)$. As it is well-known that
$Z_\gr(E(\L))=k$, in this case, it follows that
$\Ext^*_\L(\L/\rrad,\L/\rrad)$ is an infinitely generated module over
any graded subalgebra $H$ of $\HH^*(\L)$.
 
Now we introduce the first of two finite generation assumptions that
we keep throughout the paper.

\begin{assumption}[\textbf{Fg1}]
There exists a graded subalgebra $H$ of $\HH^*(\L)$ such that 
\begin{enumerate}
\item[(i)] $H$ is a commutative Noetherian ring.
\item[(ii)] $H^0=\HH^0(\L)=Z(\L)$.
\end{enumerate}
\end{assumption}
This assumption, \textbf{Fg1}, is assumed throughout the paper unless
otherwise explicitly stated.  One reason for making this assumption is
to obtain an affine variety in which to consider the support varieties
of finitely generated modules as introduced in \cite{SS}. We now make
an equivalent definition of the variety of a pair of modules $(M,N)$
in $\mod\L$ to the one given there. Let the variety of $(M,N)$ be
given by
\[V_H(M,N)=\MaxSpec(H/A_H(M,N)),\]
where $A_H(M,N)$ is the annihilator of $\Ext^*_\L(M,N)$ as an
$H$-module and this is considered as a subvariety of the homogeneous
affine variety $V_H=\MaxSpec H$. Recall that 
$V_H(M,\L/\rrad)=V_H(M,M)= V_H(\L/\rrad,M)$, and this is defined to
be the variety $V_H(M)$ of $M$. By \cite[Proposition 4.6 (b)]{SS}
$\sqrt{A_H(\L/\rrad,\L/\rrad)}$ is the ideal $\N_H$ in $H$ generated
by all homogeneous nilpotent elements in $H$. Consequently,
$V_H(\L/\rrad)=V_H$. 

To start analysing the consequences of finite generation, we need the
following proposition linking the dimension of the variety of a module
to the complexity of the module. Recall that the complexity $c_\L(M)$
of a $\L$-module $M$ is given by $c_\L(M)=\min\{b\in\mathbb{N}_0\mid
\exists a\in \mathbb{R} \text{\ such that\ } \dim_k Q_n\leq an^{b-1},
\forall n\geq 0\}$, where $\cdots \to Q_n\to Q_{n-1}\to \cdots \to
Q_0\to M\to 0$ is a minimal projective resolution of $M$. Below, 
$\gamma$ denotes the rate of growth of the dimensions of the graded
parts of a graded module.

\begin{prop}\label{prop:dim}
Let $M$ be in $\mod\L$.
\begin{enumerate}
\item[(a)] If $\Ext^*_\L(M,\L/\rrad)$ is a finitely generated
  $H$-module, then \[\dim V_H(M)=c_\L(M)<\infty.\]
\item[(b)]\sloppy If $\Ext^*_\L(\L/\rrad,M)$ is a finitely generated
  $H$-module, then \[\dim
  V_H(M)=c_\L(D(M))<\infty.\]
\end{enumerate}
\end{prop}
\begin{proof}
(a) Suppose that $\Ext^*_\L(M,\L/\rrad)$ is finitely generated as an 
$H$-module. Then using induction on the Loewy length of a module,
$\Ext_\L^*(M,N)$ is a finitely generated $H$-module for all $N$ in
$\mod\L$. In particular, $\Ext^*_\L(M,M)$ is a finitely generated
$H$-module. As the $H$-module structure on $\Ext^*_\L(M,\L/\rrad)$ 
factors through $\Ext^*_\L(M,M)$ we infer that
$\Ext^*_\L(M,\L/\rrad)$ is a 
finitely generated module over $\Ext^*_\L(M,M)$. Using these
observations the proof given in \cite[Proposition 5.3.5 and
Proposition 5.7.2]{B} carries over. In particular, the
complexity is finite. 

(b) As $V_H(M)=V_H(D(M))$ by \cite[Proposition 3.5]{SS} the claim
follows from (a).
\end{proof}

In the theory of support varieties for group rings and complete
intersections (\cite{Av,AB}) the finite generation of the extension
groups $\Ext^*(M,N)$ as modules over the ring of cohomological
operators for all modules $M$ and $N$ are of great importance. However
we first analyse finite generation for $\Ext^*_\L(M,N)$ as an 
$H$-module for specific pairs of $\L$-modules.

\begin{prop}\label{prop:localfingen}
Let $M$ be in $\mod\L$. 
\begin{enumerate}
\item[(a)] Suppose that $\Ext^*_\L(M,\L/\rrad)$ is finitely generated
as an $H$-module. If the variety of $M$ is trivial, then the projective
dimension of $M$ is finite.
\item[(b)] Suppose that $\Ext^*_\L(\L/\rrad,M)$ is finitely generated
as an $H$-module. If the variety of $M$ is trivial, then the injective
dimension of $M$ is finite.
\item[(c)] Suppose that $\Ext^*_\L(D(\L^\op),\L/\rrad)$ and
$\Ext^*_\L(\L/\rrad,\L)$ are finitely generated as $H$-modules. Then
$\L$ is a Gorenstein algebra.
\end{enumerate}
\end{prop}
\begin{proof}
(a) Since $\dim V_H(M)=c_\L(M)=0$ by Proposition \ref{prop:dim}, it is
immediate that $\Ext^n_\L(M,\L/\rrad)=(0)$
for $n\gg 0$. Therefore the projective dimension of $M$ is
finite. 

The statements in (b) and (c) follow directly from (a).
\end{proof}

This has the following immediate consequence concerning the
Generalised Naka\-yama Conjecture.

\begin{cor}
Suppose that $\Ext^i_\L(M,M\oplus \L)=(0)$ for all $i>0$.  If
$\Ext^*_\L(M,\L/\rrad)$ is finitely generated as an $H$-module, then
$M$ is projective.
\end{cor}

Any finitely generated module over a finite dimensional algebra can be
filtered by a finite filtration of semisimple modules. Next we show
that the finite generation of $\Ext^*_\L(M,N)$ for all pairs of
modules $M$ and $N$ in $\mod\L$ over a subalgebra $H$ of $\HH^*(\L)$
is equivalent to the finite generation of
$\Ext^*_\L(\L/\rrad,\L/\rrad)$ using the filtration in semisimple
modules.

\begin{prop}\label{prop:globalfingen}
The following are equivalent.
\begin{enumerate}
\item[(i)] $\Ext^*_\L(\L/\rrad,\L/\rrad)$ is finitely generated as an 
$H$-module. 
\item[(ii)] $\Ext^*_\L(M,N)$ is finitely generated as an 
$H$-module for all pairs of $\L$-modules $M$ and
$N$ in $\mod\L$.
\item[(iii)] $\HH^*(\L,B)$ is finitely generated as an $H$-module for
  all $B$ in $\mod\L^e$. 
\end{enumerate}
\end{prop}
\begin{proof}
Since $\Ext^*_\L(M,N)\simeq \HH^*(\L,\Hom_k(M,N))$, it is clear that
(iii) implies (ii) and (ii) implies (i). So it remains to
prove that (i) implies (iii). 

Assume that $\Ext^*_\L(\L/\rrad,\L/\rrad)$ is a finitely generated
$H$-module. This is equivalent to $\HH^*(\L,\Hom_k(\L/\rrad,\L/\rrad))$
being a finitely generated $H$-module. Since any simple $\L^e$-module
is isomorphic to $\Hom_k(S,T)$ for some simple $\L$-modules $S$ and
$T$ and any finitely generated $\L^e$-module is filtered in simple
$\L^e$-modules, $\HH^*(\L,B)$ is a finitely generated $H$-module for
all $B$ in $\mod\L^e$. This completes the proof.
\end{proof}
This motivates the second of our two finite generation assumptions.

\begin{assumption}[\textbf{Fg2}]
$\Ext^*_\L(\L/\rrad,\L/\rrad)$ is a finitely generated $H$-module. 
\end{assumption}
\begin{remark}
In particular note that these two assumptions \textbf{Fg1} and
\textbf{Fg2} imply that $\HH^*(\L)$ is a finitely generated
$H$-module, and consequently $\HH^*(\L)$ itself is finitely generated 
as a $k$-algebra. Similarly, $\Ext^*_\L(\L/\rrad,\L/\rrad)$ is a
finitely generated $k$-algebra. 
\end{remark}

Combining the previous result with our earlier observations in this
section we obtain the following.

\begin{thm}\label{thm:main1}
Suppose that $\L$ and $H$ satisfy \emph{\textbf{Fg1}} and
\emph{\textbf{Fg2}}. 
\begin{enumerate}
\item[(a)] The algebra $\L$ is Gorenstein.
\item[(b)] The following are equivalent for a module $M$ in $\mod\L$.
\begin{enumerate}
\item[(i)] The variety of $M$ is trivial.
\item[(ii)] The projective dimension of $M$ is finite.
\item[(iii)] The injective dimension of $M$ is finite. 
\end{enumerate}
\item[(c)] $\dim V_H(M)=c_\L(M)$ for any module $M$ in $\mod\L$. 
\end{enumerate}
\end{thm}

\section{The annihilator of $\Ext^*_\L(M,M)$}

In contrast to the previous section, the results in this section do
not need any finite generation assumptions.

From \cite{SS} the variety of a module $M$ can be defined to be
$V_H(M,M)$ for a graded subalgebra $H$ of $\HH^*(\L)$. This has the
advantage that it is given by the annihilator of $\Ext^*_\L(M,M)$,
which is the target of the graded ring homomorphism from $\HH^*(\L)$
induced by the functor $-\otimes_\L M$. Next we describe the exact
sequence of bimodules whereby the annihilator of $\Ext^*_\L(M,M)$ is
characterized. To this end let $\cdots \to P^n\to P^{n-1}\to \cdots
\to P^1\to P^0\to \L\to 0$ be a minimal projective resolution of $\L$
as a $\L^e$-module.

\begin{defin} 
Given a homogeneous element $\eta$ in $\HH^*(\L)$ of degree $n$, 
represented by a map $\eta\colon \Omega^n_{\L^e}(\L)\to\L$, we define 
the $\L^e$-module $M_{\eta}$ by the following pushout diagram
\[\xymatrix{
0 \ar[r] & \Omega^n_{\L^e}(\L) \ar[r]\ar[d]^\eta
     & P^{n-1}\ar[r]\ar[d]
     & \Omega^{n-1}_{\L^e}(\L) \ar[r]\ar@{=}[d]
     & 0\\
0 \ar[r] & \L \ar[r]^{\alpha_\eta}
     & M_\eta\ar[r]^{\beta_\eta}
     & \Omega^{n-1}_{\L^e}(\L) \ar[r]
     & 0}\]
where we denote by $\E_{\eta}$ the bottom row short exact sequence.
\end{defin} 
Note that the isomorphism class of the module $M_\eta$ is independent
of the choice of the representation of $\eta$ as a map
$\Omega_{\L^e}^n(\L)\to \L$. For a finite group $G$ the reader should
also observe that $M_\eta$ corresponds to
$\Omega_{(kG)^e}^{-1}(\Ind_{\Delta G}^{G\times G}(L_\zeta))$ for
$\zeta$ a homogeneous element in the group cohomology ring of $G$,
where $L_\zeta \simeq \Ker (\Omega^n_{kG}(k)\extto{\zeta} k)$.
Furthermore, note that $M_\eta$ is projective as a left and as a right
$\L$-module, since the same is true for $\L$ and $\Omega_{\L^e}^i(\L)$
for all $i\geq 0$. Pushing this analogy further we define $L_\eta$ to
be the module $M_\eta\otimes_\L \L/\rrad$.

Using the sequence introduced above, the elements in
$A_{\HH^*(\L)}(M,M)$ have the following characterization. We leave the
proof to the reader as it is similar to the corresponding result for
group rings (see \cite[Proposition 5.9.5]{B}).

\begin{prop}\label{prop:annihilatorelementchar}
Let $\eta$ be a homogeneous element of degree $n$ in $\HH^*(\L)$, and
let $M$ be in $\mod\L$. Then the following are equivalent.
\begin{enumerate}
\item[(i)] $\eta$ is in $A_{\HH^*(\L)}(M,M)$.
\item[(ii)] $\E_\eta\otimes_\L M$ is a split short exact sequence.
\item[(iii)] $M_\eta\otimes_\L M\simeq M\oplus\Omega_\L^{n-1}(M)\oplus
Q$ for some projective $\L$-module $Q$.
\end{enumerate}
\end{prop}

Let $\{\eta_1,\ldots,\eta_t\}$ be a finite set of homogeneous elements 
of $\HH^*(\L)$. Define the map $\alpha_{\eta_1,\ldots,\eta_t}$ to be
the map 
\[\L\otimes_\L\L\otimes_\L\cdots\otimes_\L\L\extto{\alpha_{\eta_1} 
\otimes\alpha_{\eta_2}\otimes\cdots\otimes\alpha_{\eta_t}}
M_{\eta_1}\otimes_\L M_{\eta_2}\otimes_\L\cdots\otimes_\L M_{\eta_t}.\] 
It is easy to see that this map is a monomorphism and therefore
induces an exact sequence $\E_{\eta_1,\ldots,\eta_t}$ given by 
\[0\to \L\extto{\alpha_{\eta_1,\ldots,\eta_t}} M_{\eta_1}\otimes_\L
\cdots\otimes_\L M_{\eta_t}\to X_{\eta_1,\ldots,\eta_t}\to 0.\] This
construction enables us to give a criterion for when the ideal
generated by homogeneous elements $\{\eta_1,\ldots,\eta_t\}$ is in the
annihilator $A_{\HH^*(\L)}(M,M)$.

\begin{thm}\label{thm:bimoduleofelements}
Let $\{\eta_1,\ldots,\eta_t\}$ be a finite set of homogeneous elements
in $\HH^*(\L)$, and let $M$ be in $\mod\L$.
\begin{enumerate} 
\item[(a)] The following are equivalent.
\begin{enumerate}
\item[(i)] The ideal generated by $\{\eta_1,\ldots,\eta_t\}$ is
contained in $A_{\HH^*(\L)}(M,M)$.
\item[(ii)] $\E_{\eta_i}\otimes_\L M$ is a split short exact sequence
for all $i=1,2,\ldots,t$.
\item[(iii)] $\E_{\eta_1,\ldots,\eta_t}\otimes_\L M$ is a split short
exact sequence.
\end{enumerate}
\item[(b)] \sloppy If a finite set $\{\eta_1,\ldots,\eta_t\}$ of homogeneous 
elements in $\HH^*(\L)$ is in $A_{\HH^*(\L)}(M,M)$,
then $M$ is a direct summand of $ M_{\eta_1}\otimes_\L
\cdots\otimes_\L M_{\eta_t}\otimes_\L M$.
\end{enumerate}
\end{thm}
\begin{proof} (a) By Proposition \ref{prop:annihilatorelementchar} it
remains to prove that (ii) and (iii) are equivalent. 

Assume that the exact sequence $\E_{\eta_1,\ldots,\eta_t}\otimes_\L M$
splits. Since $\alpha_{\eta_1,\ldots,\eta_t}$ can be viewed as a
composition of $\alpha_{\eta_2,\ldots,\eta_t}$ and $\alpha_{\eta_1}$,
it follows that $\E_{\eta_1}\otimes_\L M$ splits. Similarly we show
that $\E_{\eta_i}\otimes_\L M$ splits for all $i=1,2,\ldots,t$.

Assume that the exact sequence $\E_{\eta_i}\otimes_\L M$ splits for
all $i=1,2,\ldots,t$. Since $\alpha_{\eta_1,\ldots,\eta_t}$ can be
viewed as the composition of the maps
$\id_{M_{\eta_1}}\otimes_\L\cdots\otimes_\L \id_{M_{\eta_{t-1}}}\otimes
\alpha_{\eta_t}$,\ldots, $\id_{M_{\eta_1}}\otimes\alpha_{\eta_2}$ and
$\alpha_{\eta_1}$, the map $\alpha_{\eta_1,\ldots,\eta_t}\otimes\id_M$
is a composition of $t$ split monomorphisms and therefore it is a
split monomorphism itself. This proves (a).

(b) This is a direct consequence of (a).
\end{proof}

When $\L$ is selfinjective we have a further characterization of when
homogeneous elements are in $A_{\HH^*(\L)}(M,M)$. This is an immediate
consequence of Proposition \ref{prop:annihilatorelementchar} and the
fact that $\L$ is selfinjective.

\begin{prop}
Let $\L$ be a selfinjective algebra. Let $\eta$ be a homogeneous
element of $\HH^*(\L)$ of degree $n$, and let $M$ be in $\mod\L$. Then 
$\eta$ is in $A_{\HH^*(\L)}(M,M)$ if and only if $\eta\otimes_\L 
\id_M\colon \Omega^n_{\L^e}(\L)\otimes_\L M\to \L\otimes_\L M$ is zero
in $\umod\L$.
\end{prop}

We end with a result of which we saw the first glimpses in
Proposition \ref{prop:annihilatorelementchar}.

\begin{lem}
Let $\eta_1,\eta_2,\ldots,\eta_t$ be homogeneous elements in
$A_{\HH^*(\L)}(M)$. Then $M_{\eta_1}\otimes_\L \cdots \otimes_\L
M_{\eta_t}\otimes_\L M$ is in $\add\{\Omega_\L^i(M)\}_{i=0}^N\cup 
\add\L$ for some integer $N$.
\end{lem}
\begin{proof}
If $t=1$, this is the statement of Proposition
\ref{prop:annihilatorelementchar}. Suppose that $t>1$. Let
$\widetilde{M} = M_{\eta_1}\otimes_\L \cdots\otimes_\L
M_{\eta_{t-1}}\otimes_\L M$. By Proposition \ref{prop:annihilatorelementchar} 
we have $M_{\eta_t}\otimes_\L M\simeq M\oplus \Omega_\L^{\deg \eta_t-1}(M)
\oplus Q$ for some projective $\L$-module $Q$. Tensoring this
isomorphism with $M_{\eta_1}\otimes_\L\cdots\otimes_\L M_{\eta_{t-1}}$
we obtain 
\[M_{\eta_1}\otimes_\L\cdots\otimes_\L M_{\eta_t}\otimes_\L M \simeq 
\widetilde{M}\oplus M_{\eta_1}\otimes_\L\cdots\otimes_\L
M_{\eta_{t-1}}\otimes_\L \Omega_\L^{\deg \eta_t-1}(M)\oplus Q'\]
for some projective $\L$-module $Q'$. Since the middle direct summand
of the right hand side is isomorphic to
$\Omega_\L^{\deg\eta_t-1}(\widetilde{M})$ modulo projectives, the
claim now follows by induction.
\end{proof}

\section{Modules with given varieties}

In this section we return to the setting suggested by the first
section and require throughout that $\L$ satisfies \textbf{Fg1} and
\textbf{Fg2} for some graded subalgebra $H$ of $\HH^*(\L)$.

In general the variety of a module is a closed homogeneous variety. 
Here we show that any closed homogeneous variety occurs as the variety
of some module. The module we construct is not necessarily
indecomposable. 

To prove our results we make use of the bimodules $M_\eta$ introduced
in the previous section. We start by considering, for a homogeneous
element $\eta$ of positive degree in $H$, the variety of
$M_\eta\otimes_\L M$. 
\begin{prop}\label{prop:varbimodinc}
Let $\eta$ be a homogeneous element of positive degree in $H$, and let
$M$ be in $\mod\L$. 
\begin{enumerate}
\item[(a)] $V_H(M_\eta\otimes M)\subseteq V_H(L_\eta)\cap V_H(M)$.
\item[(b)] The element $\eta^2$ is in $A_H(M_\eta\otimes_\L
M,\L/\rrad)$. In particular, $V_H(L_\eta)$ is contained in 
$V_H(\langle\eta\rangle)$, and consequently $V_H(M_\eta\otimes_\L M)$ 
is contained in $V_H(\langle\eta\rangle)\cap V_H(M)$.
\item[(c)] Let $\{ \eta_1,\ldots,\eta_t\}$ be homogeneous elements in
$\HH^*(\L)$. Then $V_H(M_{\eta_1}\otimes_\L\cdots\otimes_\L
M_{\eta_t}\otimes_\L M)$ is contained in $V_H(\langle
\eta_1,\ldots,\eta_t\rangle)\cap V_H(M)$.
\end{enumerate}
\end{prop}
\begin{proof}
(a) The sequence $\E_\eta\otimes_\L M\colon 0\to M\to M_\eta\otimes_\L
M\to \Omega^{n-1}_{\L^e}(\L)\otimes_\L M\to 0$ is exact. By
\cite[Proposition 3.4]{SS} it follows that $V_H(M_\eta\otimes_\L
M)\subseteq V_H(M)$, since $\Omega^{n-1}_{\L^e}(\L)\otimes_\L M \simeq
\Omega^{n-1}_\L(M)\oplus F$ for some projective $\L$-module $F$ and
the variety is invariant under taking syzygies. 

Since the module $M$ has a finite filtration in semisimple modules,
$V_H(M_\eta\otimes_\L M)$ is contained in $V_H(L_\eta)$ by iterated
use of \cite[Proposition 3.4]{SS}. Hence the claim follows. 

(b) The proof is very similar to the group ring case. Since
\[\Ext^i_\L(\Omega_{\L^e}^{n-1}(\L)\otimes_\L M,\L/\rrad)\] is isomorphic to 
$\Ext^{i+n-1}_\L(M,\L/\rrad)$ for $i\geq 1$, the short exact sequence
\[\E_\eta\otimes_\L M\colon 0\to M\to M_\eta\otimes_\L M\to
\Omega_{\L^e}^n(\L)\otimes_\L M\to 0\] gives rise to a long exact
sequence
\begin{multline}
\cdots \to 
\Ext^i_\L(M,\L/\rrad)\extto{\cdot\eta}
\Ext^{i+n}_\L(M,\L/\rrad)\to 
\Ext^{i+1}_\L(M_\eta\otimes_\L M,\L/\rrad) \to\notag \\ 
\Ext^{i+1}_\L(M,\L/\rrad)\extto{\cdot\eta}
\Ext^{i+n+1}_\L(M,\L/\rrad)\to \cdots\end{multline}
and 
\begin{multline} 0\to \Hom_\L(\Omega_{\L^e}^{n-1}(\L)\otimes_\L M,\L/\rrad)\to 
\Hom_\L(M_\eta\otimes_\L M,\L/\rrad)\to \notag \\
\Hom_\L(M,\L/\rrad)\extto{\cdot\eta}\Ext^n_\L(M,\L/\rrad).
\end{multline}
This induces the following short exact sequences
\begin{multline} 
0\to \Hom_\L(\Omega_{\L^e}^{n-1}(\L)\otimes_\L M,\L/\rrad)\to 
\Hom_\L(M_\eta\otimes_\L M,\L/\rrad)\extto{\nu_0} \notag \\
\Ker(\cdot\eta)|_{\Hom_\L(M,\L/\rrad)}\to 0
\end{multline}
and
\begin{multline} 
0\to
\Ext^{*+n-1}_\L(M,\L/\rrad)/(\eta\Ext^{*-1}_\L(M,\L/\rrad))\extto{\mu_*}
\notag\\ \Ext^*_\L(M_\eta\otimes_\L M,\L/\rrad)\extto{\nu_*}
\Ker(\cdot\eta)|_{\Ext^*_\L(M,\L/\rrad)}\to 0\end{multline} for the
index $*$ ranging over natural numbers greater or equal to $1$. Let
$\theta$ be in $\Ext^i_\L(M_\eta\otimes_\L M,\L/\rrad)$. Then
$\nu_{i+n}(\eta\theta)=\eta\nu_i(\theta)=0$, hence $\eta\theta$ is in
$\Ker\nu_{n+i}=\Im\mu_{n+i}$. Since $\Im\mu_{n+i}$ is annihilated by
$\eta$, it follows that $\eta^2\theta=0$. The claim follows from this.

(c) This follows immediately from (b).
\end{proof}

It follows that given homogeneous elements $\{\eta_1,\ldots,\eta_t\}$
in $H$, the module $M_{\eta_1}\otimes_\L\cdots\otimes_\L
M_{\eta_t}\otimes_\L\L/\rrad$ has variety contained in
$V_H(\langle\eta_1,\ldots,\eta_t\rangle)$.  Note that the previous
result is true in general. However to show that the inclusion actually
is an equality we make full use of our assumptions \textbf{Fg1} and
\textbf{Fg2}.

We stress that throughout this section we assume the conditions
\textbf{Fg1} and \textbf{Fg2}. Recall in particular from Theorem
\ref{thm:main1} that in this case $\L$ is a Gorenstein ring.
Furthermore, these assumptions are satisfied for any block of a group
ring of a finite group (see \cite{E,V}) and more generally for a
finite dimensional cocommutative Hopf algebra (see \cite{FS}). In
addition they hold true for local finite dimensional algebras which
are complete intersections (see \cite{G}).

When $\L$ is Gorenstein the injective dimensions of $\L$ as a left and
a right module over itself are finite and they are equal, say equal to
$n$. Denote by $^\perp \L$ the full subcategory $\{ X\in \mod\L\mid
\Ext^i_\L(X,\L)=(0) \text{\ for all\ } i> 0\}$ of $\mod\L$.  Since the
variety is invariant under taking syzygies, all the different
varieties of modules occur for a module in $^\perp \L$. 

Given a module $M$ in $^\perp \L$, there is a complete resolution of
$M$
\[ \mathbb{P}_*\colon \cdots \to P_2\extto{d_2} P_1\extto{d_1}
P_0\extto{d_0} P_{-1} \extto{d_{-1}} P_{-2} \to \cdots\] where $\Im
d_0\simeq M$. This uses that $\L$ is a cotilting module; for further
details see \cite{AR}.  Define $\cExt^i_\L(M,N)$ as the homology of
$\Hom_\L(\mathbb{P}_*,N)$ at stage $i$, i.e.\ $\cExt^i_\L(M,N)=\Ker
d^*_{i+1}/\Im d^*_i$ for $i$ in $\Z$. This gives rise to the exact
sequence
\[ 0\to \P(M,N)\to \Ext^*_\L(M,N)\to \cExt^*_\L(M,N)\to
\cExt^-_\L(M,N)\to 0\] of $H$-modules, where the end terms are
$\m_\gr$-torsion modules. An $H$-module $X$ is $\m_\gr$-torsion if any
element $x$ in $X$ is annihilated by some power of $\m_\gr$.  Here
$\P(M,N)$ denotes all the homomorphisms from $M$ to $N$ factoring
through a projective module, and $\cExt^-_\L(M,N)=\oplus_{i\leq
  -1}\cExt^i_\L(M,N)$. The last fact is heavily used in the proof of
the next result, which we leave to the reader.

\begin{lem}\label{lem:localisationext}
For any maximal ideal $\frakp$ in $\MaxSpec H$ with $\frakp\neq
\m_\gr$,
\[\Ext^*_\L(M,N)_\frakp \simeq \cExt^*_\L(M,N)_\frakp\] for all
modules $M$ in $^\perp \L$ and $N$ in $\mod\L$.
\end{lem}

Next we show that $V_H(M_\eta\otimes_\L M)=V_H(\langle
\eta\rangle)\cap V_H(M)$ for any homogeneous element
$\eta$ of positive degree in $H$ and any $\L$-module $M$ in $\mod\L$. 

\begin{prop}\label{prop:varbimodequal}
Let $\eta$ be a homogeneous element of positive degree in $H$. Then
\[V_H(M_\eta\otimes_\L M)=V_H(\langle \eta\rangle)\cap 
V_H(M).\] 
In particular, $V_H(L_\eta)=V_H(\langle\eta\rangle)$. 
\end{prop}
\begin{proof}
If $V_H(\langle\eta\rangle)\cap V_H(M)$ is trivial, there is nothing
to prove. So, assume that $V_H(\langle\eta\rangle)\cap V_H(M)$ is
non-trivial. Let $n>0$ be the degree of $\eta$. 

When $\L$ is Gorenstein the syzygies $\Omega_\L^m(X)$ for a
$\L$-module $X$ are in $^\perp \L$ for $m\geq \id \L$. Since the
varieties of $M_\eta\otimes_\L M$ and $M_\eta\otimes_\L
\Omega^m_\L(M)$ are the same, for the latter is a syzygy of the
former, we can assume that $M$ is in $^\perp \L$ and therefore
$M_\eta\otimes_\L M$ is also in $^\perp \L$.

The element $\eta$ gives rise to the exact sequence
\[ 0\to \L\to M_\eta\to \Omega^{n-1}_{\L^e}(\L)\to 0\]
of $\L^e$-modules. Tensoring this sequence with $M$ and applying
$\Hom_\L(-,\L/\rrad)$ induces the long exact sequence
\begin{multline} 
\cdots\to\cExt^{i-1}_\L(M,\L/\rrad)\extto{\partial}
\cExt^i_\L(\Omega^{n-1}_{\L^e}(\L)\otimes_\L M,\L/\rrad)\to \notag\\
\cExt^i_\L(M_\eta\otimes_\L M,\L/\rrad)\to
\cExt^i_\L(M,\L/\rrad)\to\cdots\notag
\end{multline}
It is easy to see that $\partial$ is multiplication by $\eta$ when
$\cExt^i_\L(\Omega^{n-1}_{\L^e}(\L)\otimes_\L M,\L/\rrad)$ is
identified with $\cExt^{n+i}_\L(M,\L/\rrad)$. This
yields the exact sequence
\begin{multline}
\widehat{\theta}\colon 0\to \cExt^{*+n-1}_\L(M,\L/\rrad)/ 
\eta\cExt^{*-1}_\L(M,\L/\rrad)\to \notag\\
\cExt^*_\L(M_\eta\otimes_\L M,\L/\rrad)\to 
\Ker(\ \cdot\eta)|_{\cExt^*_\L(M,\L/\rrad)} \to 0.
\end{multline} 
\sloppy Choose a maximal ideal $\frakp\neq \m_\gr$ lying over $\langle
\eta, \Ann_H\Ext^*_\L(M,\L/\rrad)\rangle$. Suppose that
$\Ann_H\Ext^*_\L(M_\eta\otimes_\L M,\L/\rrad)$ is not contained in
$\frakp$. Then $\Ext^*_\L(M_\eta\otimes_\L M,\L/\rrad)_\frakp=(0)$ and
by Lemma \ref{lem:localisationext} the localisation
$\cExt^*_\L(M_\eta\otimes_\L M,\L/\rrad)_\frakp=(0)$. From the exact
sequence $\widehat{\theta}$ we infer that
$\cExt^*_\L(M,\L/\rrad)_\frakp \simeq
\eta\cExt^*_\L(M,\L/\rrad)_\frakp$. Since
$\Ext^*_\L(M,\L/\rrad)_\frakp \simeq \cExt^*_\L(M,\L/\rrad)_\frakp$ is
a finitely generated $H_\frakp$-module and $\eta$ is in $\frakp
H_\frakp$, the Nakayama Lemma implies that
$\cExt^*_\L(M,\L/\rrad)_\frakp=(0)$. Using Lemma
\ref{lem:localisationext} again $\Ext^*_\L(M,\L/\rrad)_\frakp=(0)$,
and since $\Ext^*_\L(M,\L/\rrad)$ is a finitely generated $H$-module,
the annihilator $\Ann_H\Ext^*_\L(M,\L/\rrad)$ is not contained in
$\frakp$. This is a contradiction to the choice of $\frakp$, and hence
$\Ann_H\Ext^*_\L(M_\eta\otimes_\L M,\L/\rrad)$ is contained in
$\frakp$. It follows that $V_H(\langle \eta\rangle)\cap
V_H(M)\subseteq V_H(M_\eta\otimes_\L M)$. The opposite inclusion is
proved in Proposition \ref{prop:varbimodinc}, and this completes the
proof of the proposition.
\end{proof}
The corresponding proof for a group ring of a finite group uses rank
varieties and reduction to elementary abelian subgroups, so the above
proof also gives an alternative proof in that case.

Using the above result it is easy to show that any homogeneous variety
occurs as a variety of a module. 

\begin{thm}\label{thm:allhomvar}
Let $\fraka$ be any homogeneous ideal in $H$. Then there
exists a module $M$ in $\mod\L$ such that $V_H(M)=V_H(\fraka)$.
\end{thm}
\begin{proof}\sloppy 
Suppose that $\fraka=\langle \eta_1,\eta_2,\ldots,\eta_t\rangle$ for
some homogeneous elements $\{\eta_1,\eta_2,\ldots,\eta_t\}$ in $H$. Then
$V_H(M_{\eta_1}\otimes_\L\cdots \otimes_\L
M_{\eta_t}\otimes_\L\L/\rrad)=V_H(\langle
\eta_1,\eta_2,\ldots,\eta_t\rangle)$ by Proposition
\ref{prop:varbimodequal}.
\end{proof}
Note that the module constructed for the given closed homogeneous
variety need not be indecomposable.

\section{Periodic modules}

A group ring of a finite group over a field is a symmetric algebra, so
that the Auslander-Reiten translate $\tau$ is isomorphic to
$\Omega^2$. Having a good supply of $\tau$-periodic modules gives 
information about the shape of the stable Auslander-Reiten quiver in
this case. As $\tau$-periodic and $\Omega$-periodic modules coincide
here and the $\Omega$-periodic modules are known to be controlled by
the support varieties, the theory of support varieties can be used to
(re)prove Webb's theorem (see \cite{W}). 

In this section $\L$ is a selfinjective algebra.  We take a closer
look at the construction and characterization of periodic modules. In
particular we show a generalisation of Webb's theorem for a finite
dimensional selfinjective algebra with a Nakayama functor which is of
finite order for each indecomposable module. By a periodic module we
mean throughout an $\Omega$-periodic module.

For group rings a module is periodic if and only if the variety is a
line. In fact the proof in our setting is the same as in this case, so
that we leave the details to the reader. Recall the following result
from \cite{C2} (see \cite[Proposition 5.10.2]{B}). 
\begin{prop}
Suppose $M$ is an indecomposable periodic module in $\mod\L$, of period
$n$. Then the set of nilpotent elements in
$\Ext_{\L}^{*}(M,M)$ forms an ideal, denoted by $\N_{M}$. Moreover as
vector spaces we have
\[\Ext_{\L}^{*}(M,M)\cong\N_{M}\oplus k[x]\]
where $x$ is a non-nilpotent element of degree $n$.
\end{prop}
As in the group ring case one has the following consequence, observing
that the assumptions given are sufficient (see \cite[Proposition
5.10.2]{B}).  
\begin{prop}\label{prop:periodicline}
Suppose that $\L$ and $H$ satisfy \emph{\textbf{Fg1}}. Assume that 
$M$ is an indecomposable periodic module in $\mod\L$. If $\Ext^*_\L(M,M)$ 
is a finitely generated $H$-module, then the variety of $M$ is a line.
\end{prop}

Again, in an analogous way to the group ring case, one can show that
if the variety of a module is a line, then the module is a direct sum
of periodic modules and a projective module (see \cite[Theorem
5.10.4]{B} and \cite{Ei}). For this we need to assume that $\L$ and
$H$ satisfy \textbf{Fg1} and \textbf{Fg2}, and we keep these
assumptions for the rest of this section.
\begin{thm}\label{thm:lineperiodic}
If the variety of a module $M$ in $\mod\L$ is a line, or equivalently,
if the complexity of $M$ is $1$, then $M$ is a direct sum of periodic 
modules and a projective module.
\end{thm}
Following the treatment in the group case, a homogeneous
element $\eta$ in $\HH^*(\L)$ is said to \emph{generate the 
periodicity} for a periodic module $M$ if
$V_H(\langle\eta\rangle)\cap V_H(M)$ is trivial. Thus we have yet
another analogous result (see \cite[Corollary 5.10.6]{B}).
\begin{prop}\sloppy 
If $H$ is generated as a subalgebra of $\HH^*(\L)$ by elements
$\eta_1,\eta_2, \ldots,\eta_t$ in degrees $n_1$, $n_2$,\ldots, $n_t$
and $M$ is an indecomposable periodic module in $\mod\L$, then the
period of $M$ divides one of the $n_i$.
\end{prop}

The above results can be used to construct periodic module(s) $W$ with
a variety contained in the variety of any given non-projective
indecomposable module $M$ in $\mod\L$. If $M$ is already periodic, we
can choose $W$ equal to $M$. Suppose $M$ is not periodic, or
equivalently $\dim V_H(M)> 1$. Let $\{\eta_1,\eta_2,\ldots, \eta_s\}$
be a set of homogeneous generators for $A_H(M,\L/\rrad)$. Then there
exist homogeneous elements $\{\eta_{s+1},\ldots,\eta_t\}$ in $H$ such
that the height of the ideal $\langle\eta_1,\ldots,
\eta_s,\eta_{s+1},\ldots, \eta_t\rangle$ is $\dim V_H -1$. Hence the
variety of the module $M_{\eta_1}\otimes_\L\cdots \otimes_\L
M_{\eta_t}\otimes_\L \L/\rrad$ has dimension one. By the above this is
the variety of a direct sum of periodic modules and a projective
module, so that we can, for example, choose a non-projective
indecomposable direct summand of $M_{\eta_1}\otimes_\L\cdots
\otimes_\L M_{\eta_t}\otimes_\L \L/\rrad$ to be $W$. Next we use this
to prove a generalisation of Webb's theorem.

In \cite{EH} a finite dimensional symmetric algebra $\L$ is said to
have \emph{enough} periodic modules if for any nonzero module $M$ in
$\umod\L$ there is an $\Omega$-periodic module $V$ such that
$\uHom_\L(V,M)\neq (0)$. Recall that $\uHom_\L(C,A)\simeq
D\Ext^1_\L(A,\tau^{-1}C)$ for any artin algebra and that $\tau^{-1}\simeq
\Omega^{-2}$ for a finite dimensional symmetric algebra. So, for a
symmetric algebra, showing that $\uHom_\L(V,M)$ is nonzero and showing
that $\Ext^1_\L(M,V')$ is nonzero for some periodic modules $V$ and
$V'$ are equivalent. Using support varieties we show next that finite
dimensional selfinjective algebras (satisfying \textbf{Fg1} and
\textbf{Fg2}) have enough periodic modules.

\begin{thm}
\begin{enumerate}
\item[(a)] Let $M$ be a non-projective indecomposable $\L$-module in
$\mod\L$. Then there exists a periodic $\L$-module $W$ in $\mod\L$ such that 
$\Ext^1_\L(W,M)\neq (0)$.
\item[(b)] Let $V$ be a closed homogeneous variety with
$\dim V \geq 2$. Then there exists an $\Omega$-periodic $\L$-module $W$
such that $\Ext^1_\L(W,M)\neq (0)$ for all indecomposable $\L$-modules
$M$ in $\mod\L$ with $V_H(M)=V$.
\end{enumerate}
\end{thm}
\begin{proof}
(a) If $M$ is a periodic module, we are done by choosing
$\Omega_\L^{t-1}(M)$ where $t$ is the $\Omega$-period of $M$. If $\dim
H=1$, then all indecomposable modules are periodic or projective so we 
are also done. Hence we can assume that $\dim H\geq 2$. 

Choose a homogeneous ideal $\fraka$ of height $\dim H - 1$ lying over
$A_H(M,M)$. Let $\{\eta_1,\ldots,\eta_t\}$ be homogeneous generators for
the ideal $\fraka$. Let $W=(M_{\eta_1}\otimes_\L\cdots\otimes_\L
M_{\eta_t})^*\otimes_\L \L/\rrad$. Recall that for all $i\geq 0$ 
\[\Ext^i_\L(B\otimes_\L A,C)\simeq \Ext^i_\L(A,\Hom_\L(B,C)) \simeq 
\Ext^i_\L(A,B^*\otimes_\L C)\] for $A$ and $C$ in $\mod\L$ and $B$ in
$\mod\L^e$, where $B$ is projective as a left and as a right
$\L$-module. Since $(M_{\eta_1}\otimes_\L\cdots\otimes_\L
M_{\eta_t})^*$ is projective as a left and as a right $\L$-module, we
have that
\[V_H(W) = V_H(W,\L/\rrad) = V_H(\L/\rrad,M_{\eta_1}\otimes_\L\cdots\otimes_\L
M_{\eta_t}\otimes_\L \L/\rrad)= V_H(\fraka)\] so that the variety of
$W$ is a line and therefore $W$ is a periodic module by Theorem
\ref{thm:lineperiodic}. Furthermore
$V_H(M_{\eta_1}\otimes_\L\cdots\otimes_\L M_{\eta_t}\otimes_\L M)$ is
also equal to $V_H(\fraka)$, so in particular
$M_{\eta_1}\otimes_\L\cdots\otimes_\L M_{\eta_t}\otimes_\L M$ is a
non-projective (non-injective) module. Therefore
\[ (0)\neq \Ext^1_\L(\L/\rrad,M_{\eta_1}\otimes_\L\cdots\otimes_\L
M_{\eta_t}\otimes_\L M)\simeq \Ext^1_\L(W,M).\]
The claim follows from this.

(b) Let $M$ be an indecomposable $\L$-module in $\mod\L$ with variety
$V$. As above, choose a homogeneous ideal $\fraka$ of height $\dim H -
1$ lying over $A_H(M,M)$. Then we observe that the module constructed
in (a) works for all the different modules $M$ with $V_H(M)=V$.
\end{proof}

This gives the promised generalisation of Webb's theorem. 

\begin{thm}
Suppose that the Nakayama functor is of finite order on any
indecomposable module in $\mod\L$. Then the tree class of a component
of the stable Auslander-Reiten quiver of $\L$ is one of the following:
a finite Dynkin diagram ($\mathbb{A}_n$, $\mathbb{D}_n$,
$\mathbb{E}_{6,7,8}$), an infinite Dynkin diagram of the type
$\mathbb{A}_\infty$, $\mathbb{D}_\infty$, $\mathbb{A}^\infty_\infty$
or a Euclidean diagram.
\end{thm}
\begin{proof}
The Auslander-Reiten translate for a finite dimensional selfinjective
algebra is the composition of the Nakayama functor $\N$ with the second
syzygy $\Omega^2$, where these two functors commute.  Hence, by
assumption an indecomposable module is $\tau$-periodic if and only if
the module is $\Omega$-periodic. We use this and the above result to
construct subadditive functions on the tree class of the stable
components of $\L$.

Let $\C$ be a component of the stable Auslander-Reiten quiver $\G_s$.
In Theorem 3.7 in \cite{SS} it is shown that all indecomposable
modules in $\C$ have the same variety, say $V$. If $\dim V=1$, then
all modules in the component are periodic. Then the tree class is a
finite Dynkin diagram or $\mathbb{A}_\infty$ by \cite{HPR}.

If $\dim V\geq 2$, then no module in the component $\C$ is periodic,
and by the previous result there exists a periodic module $W$ such
that $\Ext^1_\L(W,M)\neq (0)$ for all $M$ in $\C$. We can assume
without loss of generality that $W$ is $\tau$-periodic of period one,
that is, $W\simeq \tau W$. Define $f(M)=\dim_k\Ext^1_\L(W,M)$ for any
$M$ in $\C$. Then
\begin{align}
f(\tau M) & =\dim_k\Ext^1_\L(W,\tau M)\notag\\
          & =\dim_k\Ext^1_\L(\tau^{-1}W,M)\notag\\
          & =\dim_k\Ext^1_\L(W,M)=f(M).\notag
\end{align}
Now this gives rise to a subadditive function on the tree class of
$\C$, and it follows from \cite{HPR} that the tree class of $\C$ is
one of diagrams listed above.
\end{proof}
\begin{remark}
Note that when the dimension of the variety of a module $M$ in the
component $\C$ in the previous proof is at least two, then the
function constructed is actually additive. To see this, observe that
none of the indecomposable non-projective direct summands of $W$ can
lie in $\C$ as their variety has dimension one. If $0\to A\extto{f} 
B\extto{g} C\to 0$ is an almost split sequence in $\C$, then the sequence 
\[0\to \Ext^1_\L(W,A)\extto{f^*}\Ext^1_\L(W,B)\extto{g^*} 
\Ext^1_\L(W,C)\to 0\]
is exact as the kernel of $f^*$ and the cokernel of $g^*$ are 
isomorphic to $\delta^*(\Omega^i_\L(W))$ for $i=0$ and $i=1$,
respectively. Here $\delta^*$ is given by the exact sequence of the
functor
\[0\to \Hom_\L(-,A)\to \Hom_\L(-,B)\to \Hom_\L(-,C)\to \delta^*\to 0,\] 
and $\delta^*(X)\neq (0)$ if and only if $C$ is a direct summand of
$X$. 
\end{remark}

\section{Representation type and complexity}

Here we give a brief discussion on relationships between
representation type and complexity of modules over a selfinjective
algebra $\L$. 

It was first observed by Heller in \cite{H} that if $\L$ is of finite
representation type, then all the indecomposable non-projective
modules are $\Omega$-periodic and therefore all of complexity one. The
converse of this statement is not true, since there exist finite
dimensional preprojective algebras of wild representation type with
$\L$ being a periodic $\L^e$-module (and consequently all
indecomposable non-projective modules are periodic and of complexity
one).

If $\L$ is of tame representation type, then it is shown by Rickard in
\cite{R} that all indecomposable non-projective modules have
complexity at most two. By using the same example as above the
converse is also not true here. However a partial converse is known,
and we include a proof here for completeness. 

\begin{prop}
Suppose $\L$ and $H$ satisfy \emph{\textbf{Fg1}} and
\emph{\textbf{Fg2}} with $\dim H\geq 2$. Then $\L$ is of infinite
representation type and $\L$ has an infinite number of indecomposable
periodic modules lying in infinitely many different components of the
stable Auslander-Reiten quiver. 
\end{prop}
\begin{proof}
Suppose that $\dim H=d\geq 2$. Then by the Noether Normalisation Theorem
there exists a polynomial ring $k[x_1,\ldots,x_d]$ generated by
homogeneous elements $x_i$ in $H$ of degrees $n_1$, $n_2$,\ldots,
$n_d$, respectively, over which $H$ is a finitely generated module. 
Choose natural numbers $r$ and $s$ with $r+s$
minimal such that $x_1^r$ and $x_2^s$ have the same degree. Let
$\eta_\alpha=\alpha_1x_1^r+\alpha_2x_2^s$ for
$\alpha=(\alpha_1,\alpha_2)$ in $\mathbb{P}^1(k)$. Consider the
modules $C_\alpha=M_{\eta_\alpha}\otimes_\L
M_{x_3}\otimes_\L\cdots \otimes_\L M_{x_d}\otimes_\L \L/\rrad$. Then
$V_H(C_\alpha)=V_H(\langle \eta_\alpha,x_3,\ldots,x_d\rangle)$. It is
easy to show that $\eta_\alpha$ is an irreducible element in $H$ when
$\alpha$ is different from $(1,0)$ and $(0,1)$, so that
$V_H(C_\alpha)$ is an irreducible variety ($\sqrt{\langle
  \eta_\alpha,x_3,\ldots,x_d\rangle}$ is a prime ideal). If $X$ is any
indecomposable non-projective direct summand of $C_\alpha$, then
$V_H(X)$ is a closed subvariety of $V_H(C_\alpha)$.  We infer that
$V_H(X)=V_H(C_\alpha)$. We can then construct a
$\mathbb{P}^1(k)$-family of indecomposable periodic modules
$\{X_\alpha\}$ by choosing for each $\alpha$ in $\mathbb{P}^1(k)$ an
indecomposable non-projective direct summand $X_\alpha$ of $C_\alpha$.
Since $V_H(X_\alpha)\neq V_H(X_{\alpha'})$ for $\alpha\neq \alpha'$ in
$\mathbb{P}^1(k)$, the claim follows as the field is infinite.
\end{proof}

\section{The variety of an indecomposable module is connected}

Throughout this section we assume that $\L$ satisfies \textbf{Fg1} and
\textbf{Fg2} for some graded subalgebra $H$ of $\HH^*(\L)$. This
section is devoted to showing that the variety of an indecomposable
module is connected. Here $\id X$ denotes the injective dimension of a
module $X$. 

To show the result mentioned above we need preliminary results, and we
first give a sufficient condition for the vanishing of all high enough
extension groups between two modules.

\begin{prop}\label{prop:intersectvar}
Let $M$ and $N$ be two $\L$-modules in $\mod\L$. Suppose that
$V_H(M)\cap V_H(N)$ is trivial. Then $\Ext^i_\L(M,N)=(0)$ for all
$i > \id \L_\L$.

In particular, if $\L$ is selfinjective, then $\Ext^i_\L(M,N)=(0)$ for
all $i\geq 1$.
\end{prop}
\begin{proof}
Suppose that $A_H(M,M)$ is generated by $\{\eta_1,\ldots,\eta_t\}$.  By
Theorem \ref{thm:bimoduleofelements} the module $M$ is a direct
summand of $M_{\eta_1}\otimes_\L\cdots \otimes_\L M_{\eta_t}\otimes_\L
M$. Therefore $\Ext^i_\L(M,N)$ is a direct summand of
$\Ext^i_\L(M_{\eta_1}\otimes_\L\cdots \otimes_\L M_{\eta_t}\otimes_\L
M,N)$. We have that
\[\Ext^i_\L(M_{\eta_1}\otimes_\L\cdots \otimes_\L
M_{\eta_t}\otimes_\L M,N)\simeq
\Ext^i_\L(M, (M_{\eta_1}\otimes_\L\cdots \otimes_\L
M_{\eta_t})^*\otimes_\L N).\]
Moreover 
\[\Ext^i_\L(\L/\rrad,(M_{\eta_1}\otimes_\L\cdots \otimes_\L
M_{\eta_t})^*\otimes_\L N)\simeq
\Ext^i_\L(M_{\eta_1}\otimes_\L\cdots \otimes_\L M_{\eta_t}\otimes_\L
\L/\rrad,N)\] so that the variety of $(M_{\eta_1}\otimes_\L\cdots
\otimes_\L M_{\eta_t})^*\otimes_\L N$ is contained in the
intersection of the varieties of $M_{\eta_1}\otimes_\L\cdots
\otimes_\L M_{\eta_t}\otimes_\L \L/\rrad$ and of $N$. The variety of
$M_{\eta_1}\otimes_\L\cdots \otimes_\L M_{\eta_t}\otimes_\L \L/\rrad$
is contained in that of $M$, so that by assumption the variety of
$(M_{\eta_1}\otimes_\L\cdots \otimes_\L M_{\eta_t})^*\otimes_\L N$
is trivial. 

By Theorem \ref{thm:main1} $\L$ is Gorenstein and $\id
(M_{\eta_1}\otimes_\L\cdots \otimes_\L M_{\eta_t})^*\otimes_\L N \leq
\id \L_\L$, so that $\Ext^i_\L(M,N)=(0)$ for all $i > \id\L_\L$.

Since $\id \L_\L=0$ when $\L$ is selfinjective, the last claim is
clear. 
\end{proof}

In a similar way to the group case we show that the variety of an
indecomposable module is connected. The next result is a crucial step
in the proof and the assumptions \textbf{Fg1} and \textbf{Fg2} are not
needed for this. 
\begin{lem}\label{lem:prodseq}
Given two homogeneous elements $\eta_1$ and $\eta_2$ of positive degree
in $\HH^*(\L)$, there is an exact sequence
\[ 0\to \Omega_{\L^e}^n(M_{\eta_1}) \to M_{\eta_2\eta_1}\oplus F\to
M_{\eta_2}\to 0\]
of $\L^e$-modules, where $F$ is projective and $n$ is the degree of
$\eta_2$. 
\end{lem}
\begin{proof}
Given the map $\eta_1\colon \Omega^m_{\L^e}(\L)\to \L$ we obtain the
exact sequence 
\[\Theta\colon 0\to \Omega_{\L^e}^1(M_{\eta_1})\oplus Q'\to
\Omega_{\L^e}^m(\L)\oplus Q\extto{\left(\begin{smallmatrix}
\eta_1\\
f\end{smallmatrix}\right)} \L\to 0\]
for some projective $\L^e$-modules $Q$ and $Q'$. Using that 
$\eta_2\eta_1$ is given by the composition 
$\eta_2\Omega_{\L^e}^n(\eta_1)$, this gives rise to the exact 
commutative diagram
\[\xymatrix{
0\ar[r] & \L\ar[r] \ar@{=}[d] & M_{\eta_2\eta_1} \ar[r] \ar[d] &
\Omega_{\L^e}^{m+n-1}(\L)\ar[r]\ar[d]^{\Omega_{\L^e}^{n-1}(\eta_1)} & 0\\
0\ar[r] & \L\ar[r] & M_{\eta_2} \ar[r] & 
\Omega_{\L^e}^{n-1}(\L)\ar[r] & 0}\]
By applying the Horseshoe Lemma to the sequence $\Theta$ we know that
we can add a projective $\L^e$-module $F$ to make the map
$\Omega_{\L^e}^{m+n-1}(\L)\oplus F\extto{\left(\begin{smallmatrix}
\Omega_{\L^e}^{n-1}(\eta_1)\\
f'\end{smallmatrix}\right)} \Omega_{\L^e}^{n-1}(\L)$ onto with kernel
isomorphic to $\Omega_{\L^e}^n(\L)$. By adding $F$ also to
$M_{\eta_2\eta_1}$ and making the necessary adjustments to the maps,
we obtain the desired sequence. 
\end{proof}

By combining Proposition \ref{prop:annihilatorelementchar},
Proposition \ref{prop:intersectvar} and Lemma \ref{lem:prodseq} and
making the appropriate modifications to the proof in the group ring
case we obtain that the variety of an indecomposable is connected (see
\cite{C3}, \cite[Theorem 5.12.1]{B}).
\begin{thm}
Let $M$ be in $\mod\L$. If $V_H(M)=V_1\cup V_2$ 
for some homogeneous non-trivial varieties $V_1$ and $V_2$ with 
$V_1\cap V_2$ trivial, then $M\simeq M_1\oplus M_2$ with
$V_H(M_1)=V_1$ and $V_H(M_2)=V_2$.
\end{thm}

\end{document}